\documentclass[a4paper]{article}
\usepackage{lmodern}
\usepackage[T1]{fontenc}
\usepackage{hyperref}
\usepackage{amsfonts,amsmath,amssymb,amsopn}
\usepackage[all]{xy}
\usepackage{vmargin}
\usepackage{tikz}
\title{On quantum state of numbers}
\author{Bernard Le Stum \& Adolfo Quir\'os\footnote{Supported by grants MTM2009-07291 from Ministerio de Ciencia e Innovaci\'on (Spain) and MTM2012-35849 from Ministerio de Econom\'{\i}a y Competitividad (Spain)}}
\date{Version of \today}
\newtheorem{thm}{Theorem}[section]
\newtheorem{prop}[thm] {Proposition}
\newtheorem{cor}[thm] {Corollary}
\newtheorem{lem}[thm] {Lemma}
\newtheorem{dfn}[thm] {Definition}

\newenvironment{xmp}[1][Examples.]{\begin{trivlist} \item[\hskip \labelsep {\bfseries #1}]}{\end{trivlist}}
\newenvironment{pf}[1][Proof]{\begin{trivlist} \item[\hskip \labelsep {\bfseries #1}]}{\end{trivlist}}
\newenvironment{rmk}[1][Remarks.]{\begin{trivlist} \item[\hskip \labelsep {\bfseries #1}]}{\end{trivlist}}
\begin{document}

\maketitle

\bigskip

\begin{center}
\textbf{Abstract}
\end{center}

We introduce the notions of quantum characteristic and quantum flatness for arbitrary rings.
More generally, we develop the theory of quantum integers in a ring and show that the hypothesis of quantum flatness together with positive quantum characteristic generalizes the usual notion of prime positive characteristic.
We also explain how one can define quantum rational numbers in a ring and introduce the notion of twisted powers.
These results play an important role in many different areas of mathematics and will also be quite useful in a subsequent work of the authors.

\tableofcontents

\section*{Introduction}

Quantum mathematics is obtained by making a \emph{small} perturbation $q$ on a usual mathematical object giving rise to its \emph{$q$-analog}.
Alternatively, one may consider the full collection of objects obtained for different values of $q$, giving rise to different \emph{$q$-states} of the same usual object.
We are interested here in the perturbations of the unit of a ring.
And by \emph{small}, we mean that $q$ is a non trivial root of unity.
Actually, the same process works for any perturbations $q$, which can be \emph{big} ($q$ transcendental) or even \emph{trivial} ($q = 1$).
One should not use the word \emph{quantum} in this more general situation and instead say \emph{twisted} for example.
But we will not do this here because we are actually interested in the former case in the end.

Applying this principle to a usual number, we may consider the various quantum states of the realizations of this number.
More precisely, if $R$ is a ring with unit $1$, one defines for any $q \in R$, the $q$-states of the first natural integers as follows:
\begin{displaymath}
(0)_{q} = 0, \quad (1)_{q} = 1, \quad (2)_{q} = 1 + q, \quad (3)_{q} = 1 + q + q^2, \quad \ldots.
\end{displaymath}

When $R = \mathbb Z$ is the ring of integers and $q = 1$, we recover usual natural numbers ; but if we allow some other $q > 0$, we get  the so called $q$-integers in $\mathbb Z$.
These $q$-integers may be used to develop $q$-combinatorics.
For example the $q$-analog of a binomial coefficient will count the number of rational points of the corresponding Grassmanian over a finite field with $q$ elements (\cite{KacCheung02}, Theorem 7.1).

When $R = \mathbb Z[t]$ is the polynomial ring over the integers and $q = t$, then the $q$-analog of binomial coefficients are given by the \emph{Gaussian polynomials} (these rational functions do live inside $\mathbb Z[t]$):
\begin{displaymath}
{n \choose k}_{t} := \frac {(1 - t^n)(1 - t^{n-1}) \cdots (1 -t^{n-k+1})}{(1 - t^k)(1 - t^{k-1}) \cdots (1 - t)} \in \mathbb Q(t).
\end{displaymath}
This case is very important in the theory of integer partitions (Ramanujan generating $q$-series).
See, for example, chapter 3 of \cite{Andrews76} or section 1.8 in \cite{Stanley12}.

Note that we can consider the case $R = \mathbb Z[t]$ or $\mathbb Q(t)$, and $q = t$, as the generic situation and many formulas that will be valid for any ring $R$ and any $q$ can be recovered from this particular case.

When $R = \mathbb C$ and $q \neq 1$, we may more generally define the $q$-state of any complex number $a$ once we make the choice of a branch of the logarithm that is defined at $q$:
\begin{displaymath}
(a)_{q} = \frac {1 - \exp\left(a\log(q)\right)}{1 - q} \in \mathbb C
\end{displaymath}
(when $R = \mathbb R$ and $q > 0$, we can use the usual logarithm).
When $|q| \neq 1$, we enter the realm of $q$-difference equations (see \cite{DiVizioRamisSauloyZhang03} for example).
When $q$ is a non trivial root of unity, then we get the numbers that appear in the theory of quantum groups (see \cite{Kassel95} for example).
Actually, the subjects in which $q$-analogs are fruitful keep expanding, from $q$-hypergeometric series (see \cite{Ernst12} for a thorough treatment of $q$-Calculus or \cite{KacCheung02} for a more concise introduction) to Number Theory \cite{BorisovNathansonWang04} or even Multiple $q$-Zeta Values \cite{Bradley05}.

Note that if $R$ is a ring of characteristic $p > 0$ and $q = 1$, we will have $(p)_{q} = 0$.
Also, if $q$ is a primitive $p$-th root of unity for some integer $p \geq 2$, we will have $(p)_{q} = 0$ whatever the characteristic of $R$ is.
Therefore, it appears that from a quantum point of view, roots of unity and positive characteristic share a common property.
Starting from this consideration, one may want to lift to characteristic zero some results that are already known in characteristic $p > 0$ at the cost of replacing usual mathematical objects by their $q$-analog where $q$ is a root of unity.
Michel Gros and the first author have been successful in doing this in \cite{GrosLeStum13} but we want to investigate this relation in more details in the future.
For example, the three of us are developing a quantum confluence theory and will introduce quantum divided powers.

The purpose of this article is to present many properties of quantum numbers in a complete and general form with full proofs.
Most - if not all - formulas can be found elsewhere in the literature (and this is particularly true for the formulas of section \ref{binomial} that have been well known for a long time).
However, they are usually stated with unnecessary hypothesis and their proofs often do not extend to the general case.
We wish that our presentation will provide a quick and easy reference for the mathematical community.

In section \ref{integer}, we define quantum integers (or more precisely, the quantum states of an integer) in a ring and study how the choice of the data will affect the behavior of those quantum integers.
In particular, we introduce the notion of quantum characteristic and quantum flatness.

In section \ref{binomial}, we define quantum factorials and quantum binomial coefficients.
Then we state and prove some classical results on binomial coefficients with a special emphasis on Lucas formula: it is valid under the assumptions of finite quantum characteristic and quantum flatness.

In section \ref{rational}, we define the quantum state of a rational number.
This seems new to us.
Instead of choosing a branch of the logarithm as in the complex case, one need to make a compatible choice of roots.
 We will explain this in detail.
 
In section \ref{twisted}, we consider a commutative algebra endowed with an endomorphism and introduce the notion of twisted powers.
We show that in the case of a dilatation, we recover some of the formulas that were obtained in the previous sections.
 
We wish to thank Michel Gros with whom we had many conversations related to the notions that are developed here.

Throughout this article, $R$ denotes an associative ring with unit and $q$ is an element of $R$.

\section{Quantum integers} \label{integer}

\begin{dfn}
If $m \in \mathbb N$, the \emph{$q$-state} (also called \emph{quantum state} when $q$ is part of the data) of $m$ is
\begin{displaymath}
(m)_{q} = \sum_{i=0}^{m-1} q^i \in R.
\end{displaymath}
We will also say that $(m)_{q}$ is a \emph{$q$-integer} (or a \emph{quantum integer}) of $R$.

If $q$ is invertible in $R$ and $m \neq 0$, we define the $q$-state of $-m$ as
\begin{displaymath}
(-m)_{q} = - \sum_{i=1}^{m} q^{-i} \in R
\end{displaymath}
and we will also call $(-m)_{q}$ a \emph{$q$-integer}.
\end{dfn}

In other words, we have
\begin{displaymath}
(0)_{q} = 0, \quad (1)_{q} = 1, \quad (2)_{q} = 1+q, \quad \ldots , \quad (m)_{q} = 1 + q + \cdots + q^{m-1}, \quad \ldots
\end{displaymath}
and when $q \in R^\times$,
\begin{displaymath}
(-1)_{q} = - \frac 1q, \quad (-2)_{q} = - \frac 1q - \frac 1{q^2} = - \frac {1 +q}{q^2}, \quad \ldots  ,
\end{displaymath}
\begin{displaymath}
(-m)_{q} = - \frac 1q - \cdots - \frac 1{q^m} = - \frac {1 + q + \cdots + q^{m-1}}{q^m}, \quad \ldots
\end{displaymath}

Alternatively, one may define $(m)_{q}$ by induction on $m$ as follows:
\begin{displaymath}
(0)_{q} := 0 \quad \mathrm{and} \quad (m+1)_{q} := (m)_{q} + q^m.
\end{displaymath}
One may also define $(m)_{q} := -q^{m}(-m)_{q}$ for $m < 0$ when $q \in R^\times$.

\begin{rmk}
\begin{enumerate}
\item 
These formulas take place inside the subring of $R$ generated by $q$ (and $q^{-1}$ if $q \in R^\times$).
And this last ring is commutative.
In particular, we should not worry to much about $R$ being commutative or not.
\item When $u : R \to R'$ is a ring homomorphism with $u(q) = q'$, we have
\begin{displaymath}
\forall m \in \mathbb N, \quad (m)_{q'} = u((m)_{q})
\end{displaymath}
(and the same result for $m < 0$ when $q \in R^\times$).
Using this property, we can reduce many (but not all) proofs, first to the case $R = \mathbb Z[t]$ and $q = t$, and then even to the case $R = \mathbb Q(t)$ and $q = t$.

\item
When $v \in R^\times$, one also defines the \emph{symmetric quantum state} of $n \in \mathbb Z$ (see for example section 1.3.3 of \cite{Lusztig10} or section VI.1 of \cite{Kassel95}) by
\begin{displaymath}
[n]_{v} = \frac {v^n -v^{-n}}{v - v^{-1}}.
\end{displaymath}
One can easily check that
\begin{displaymath}
[n]_{v} = \frac {(n)_{v^2}}{v^{n-1}}.
\end{displaymath}
It follows that almost any formula from one theory can be translated into the other one.
\end{enumerate}
\end{rmk}

\begin{xmp}
\begin{enumerate}
\item For $R = \mathbb Q(t)$ and $q = t$, we have
\begin{displaymath}
(m)_{q} = \frac {1 - t^m}{1 -t}.
\end{displaymath}
\item When $q = 1_{R}$ is the unit of $R$ (we will just say $q = 1$ in the future), we have $(m)_{q} = m1_{R}$.
And the canonical map $\mathbb Z \to R$ induces a bijection
\begin{displaymath}
\xymatrix@R0cm{
\mathbb Z/p\mathbb Z \ar@{<->}[r] & \{q-\mathrm{integers}\ \mathrm{in}\ R\}
}
\end{displaymath}
where $p := \mathrm{Char}(R)$.
\item For $R = \mathbb C$ and $q = e^{\frac{2\pi\sqrt{-1}}p}$ with $p \in \mathbb N \setminus \{0\}$, we obtain again a bijection
\begin{displaymath}
\xymatrix@R0cm{
\mathbb Z/p\mathbb Z  \ar@{<->}[r] & \{q-\mathrm{integers}\ \mathrm{in}\ R\}.
}
\end{displaymath}
This is illustrated in the case $p=5$ as follows:
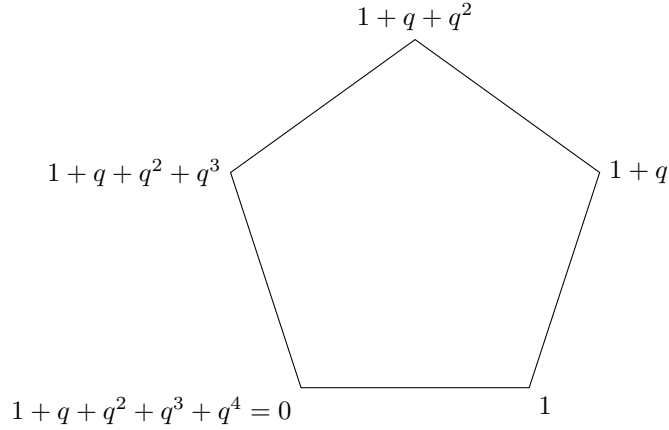
\begin{center}
\begin{tikzpicture}[scale=3]
\draw (1,0) node[below right]{$1$} -- ++ (72:1) node[right]{$1 + q$} -- ++ (2*72:1) node[above]{$1 + q + q^2$} -- ++ (3*72:1) node[left]{$1 + q + q^2+ q^3$} -- ++ (4*72:1)node[below left]{$1 + q + q^2+ q^3 + q^4 = 0$}-- cycle;
\end{tikzpicture}
\end{center}
\end{enumerate}
\end{xmp}

The following result is immediate but very important:

\begin{lem} \label{explicit}
For all $m \in \mathbb N$ (or $m \in \mathbb Z$ when $q \in R^\times$), we have
\begin{displaymath}
(1 -q)(m)_{q} = 1 - q^m.
\end{displaymath}
In particular, if $1 - q$ is invertible in $R$, we have
\begin{equation} \label{analog}
(m)_{q} = \frac {1 - q^m}{1 -q}.
\end{equation}
\end{lem}

\begin{pf}
For $m \in \mathbb N$, we have
\begin{displaymath}
(1 -q)(m)_{q} = (1 -q) \sum_{i=0}^{m-1} q^i = \sum_{i=0}^{m-1} q^i  - \sum_{i=1}^{m} q^i =1 - q^m,
\end{displaymath}
and when $q \in R^\times$,
\begin{displaymath}
(1 -q)(-m)_{q} = - (1 -q) \sum_{i=1}^{m} q^{-i} = - \sum_{i=1}^{m} q^{-i}  + \sum_{i=0}^{m-1} q^{-i} =1 - q^{-m}. \quad \Box
\end{displaymath}
\end{pf}

Note that the condition of the second assertion in the lemma implies that $q \neq 1$.
Conversely, if $q \neq 1$ and $q$ belongs to some subfield $K$ of $R$, then the condition is fulfilled.
This will often be the case in practice and formula \eqref{analog} is frequently used as an alternative definition for $q$-integers.

\begin{prop} \label{arit}
For all $m, n \in \mathbb N$ (or $\mathbb Z$ when $q \in R^\times$), we have
\begin{equation}\label{formul}
(m+n)_{q} = (m)_{q} + q^m(n)_{q}
\end{equation}
and
\begin{equation}\label{formul2}
(mn)_{q} = (m)_{q} (n)_{q^m}.
\end{equation}
\end{prop}

\begin{pf}
Pulling back along the canonical map $\mathbb Z[t] \to R$ (or $\mathbb Z[t,t^{-1}] \to R$ when $q \in R^\times$) that sends $t$ to $q$, we may first assume that $R = \mathbb Z[t]$ (or $R =\mathbb Z[t,t^{-1}]$ in the second case) and $q = t$.
Then, pushing through the embedding of $R$ into $\mathbb Q(t)$, we may actually assume that $R = \mathbb Q(t)$ (and still $q = t$).
Then, the formulas read
\begin{displaymath}
\frac {1 -t^{m+n}}{1 -t} = \frac {1 -t^m}{1 -t} + t^m \times \frac {1 -t^n}{1 -t} 
\end{displaymath}
and
\begin{displaymath}
\frac {1 -t^{mn}}{1 -t} = \frac {1 - t^m}{1 -t} \times \frac {1 - (t^m)^n}{1 -t^m} \quad \mathrm{for}\ m\neq 0.
\end{displaymath}
Of course, for $m = 0$, we have $(mn)_{q} = (0)_{q} = 0$ and also $(m)_{q} (n)_{q^m} = (0)_{q}(n)_{1} = 0 \times n = 0$.
 $\quad \Box$\end{pf}
 
\begin{dfn}
The \emph{$q$-characteristic} (or \emph{quantum characteristic} when $q$ is fixed) of $R$ is the smallest positive integer $p$ such that $(p)_{q} = 0$ if it exists and $0$ otherwise.
We will then write  $q\mathrm{-char}(R) = p$.
\end{dfn}

\begin{xmp}
\begin{enumerate}
\item Assume that $q = 1$.
Then, the quantum characteristic is the \emph{usual} characteristic of the ring $R$.
\item If $R = K[t]$ is a polynomial ring over a commutative ring $K$, and $q = t$, then $q\mathrm{-char}(R) = 0$.
\item If $R = \mathbb C$ and $q = e^{\frac {2\pi\sqrt{-1}}p}$, then the $q$-characteristic of $R$ is $p$.
\item Assume $R = \mathbb Z/n\mathbb Z$ and  $1 \neq q = \bar m \in R$.
Then the reader can check that
\begin{displaymath}
q\mathrm{-char}(R) > 0 \Leftrightarrow q \in R^\times.
\end{displaymath}
More precisely, one shows that the $q$-characteristic of $R$ is the order of $m$ in $(\mathbb Z/(m-1)n\mathbb Z)^\times$.
\end{enumerate}
\end{xmp}

\begin{prop} Let $p$ be a positive integer.
If $q\mathrm{-char}(R) = p$, then $q$ is a $p$-th a root of unity.
In particular, $q$ is invertible.
\end{prop}

\begin{pf}
We have $1 - q^p = (1 -q)(p)_{q}$.
Thus, if $(p)_{q} = 0$, we have $q^p = 1$ and $q$ is a root of unity.
$\quad \Box$
\end{pf}

\begin{prop}
If $q\mathrm{-char}(R) = p$, then the set of $m \in \mathbb N$ such that $(m)_{q} = 0$ is exactly the monoid $p\mathbb N$.
\end{prop}

If we allow $m < 0$ when $q \in R^\times$, then we get $p\mathbb Z$.

\begin{pf}
When $p = 0$, this is clear.
If $p > 0$, we can always write $m = np + r$ with $0 \leq r < p$ and $n \in \mathbb N$.
Using proposition \ref{arit}, one sees that
\begin{displaymath}
(m)_{q} = q^r(np)_{q} + (r)_{q} = q^r (n)_{q^p}(p)_{q} + (r)_{q} = (r)_{q}
\end{displaymath}
and therefore
\begin{displaymath}
(m)_{q} = 0 \Leftrightarrow (r)_{q} = 0 \Leftrightarrow r = 0 \Leftrightarrow m \in p\mathbb N. \quad \Box
\end{displaymath}
\end{pf}

\begin{prop} \label{divp}
Assume that $q\mathrm{-char}(R) = p > 0$.
Then we have
\begin{enumerate}
\item If $m, n \in \mathbb Z$ satisfy $m \equiv n \mod p$, then $(m)_{q} = (n)_{q}$.
\item If $m \in \mathbb Z$ is prime with $p$, then $(m)_{q}$ is invertible.
\end{enumerate}
\end{prop}

\begin{pf}
We will use both the fact that $(p)_{q} = 0$ and its immediate consequence $q^p = 1$.
For the first assertion, we may write $m = pv+n$ with $v \in \mathbb Z$.
We obtain
\begin{displaymath}
(m)_{q} = (pv+n)_{q} = (p)_{q}(v)_{q^p} + q^{pv}(n)_{q} = (n)_{q}.
\end{displaymath}
For the second one, we may write $mu =pv+1$ with $u, v \in \mathbb Z$ and we get
\begin{displaymath}
(m)_{q}(u)_{q^m} = (mu)_{q} = (pv+1)_{q} = (p)_{q}(v)_{q^p} + q^{pv}(1)_{q} = 1. \quad \Box
\end{displaymath}
\end{pf}

For further use, we prove the following:

\begin{lem} \label{pqzer}
Let $m \in \mathbb N \setminus \{0\}$.
Assume $R$ has no $(m)_{q}$-torsion.
Then, we have $(m)_{q} = 0$ if and only if one of the following conditions is fulfilled:
\begin{enumerate}
\item $q$ is a non trivial $m$-th root of unity.
\item $\mathrm{Char} (R) \mid m$ and $q = 1$.
\end{enumerate}
In both cases, $q$ is an $m$-th root of unity and, in particular, it is invertible.
\end{lem}

Note that the lemma is also valid for $m < 0$ when $q \in R^\times$.

\begin{pf}
Our hypothesis means that
\begin{displaymath}
(m)_{q}a = 0 \Leftrightarrow ((m)_{q} = 0\ \mathrm{or}\ a = 0).
\end{displaymath}
Since $(m)_{q}(1 - q) = 1 - q^m$, we see that $q$ is an $m$-th root of unity if and only if $(m)_{q} = 0$ or $q = 1$.
In the case $q \neq 1$, we obtain that $(m)_{q} = 0$ if and only if $q$ is an $m$-th root of unity.
When $q = 1$, the quantum state of $m$ is $m$ itself, but seen as an element of $R$.
In particular, $(m)_{1} = 0$ if and only if $p \mid m$ where $p = \mathrm{Char} (R)$.
$\quad \Box$
\end{pf}

\begin{dfn}
\begin{enumerate}
\item
The ring $R$ is said to be \emph{$q$-flat} (or \emph{quantum-flat} when the reference to $q$ is clear) if $R$ has no $(m)_{q}$-torsion for any $m \in \mathbb N$.
\item
The ring $R$ is said to be \emph{$q$-divisible} (or \emph{quantum-divisible} when the reference to $q$ is clear) if $(m)_{q} \in R^\times$ whenever $(m)_{q} \neq 0$.
\end{enumerate}
\end{dfn}

Saying that $R$ is $q$-flat means that
\begin{displaymath}
\forall m \in \mathbb N, a \in R, \quad (m)_{q}a= 0 \Leftrightarrow ((m)_{q} = 0\ \mathrm{or}\ a = 0).
\end{displaymath}
Note that the condition will then also hold for $m < 0$ when $q \in R^\times$.

Of course, $q$-divisibility always implies $q$-flatness.

\begin{xmp}
\begin{enumerate}
\item If $R$ is an integral domain (resp. a field), it is $q$-flat (resp. $q$-divisible) whatever $q$ is.
In particular, if $R = \mathbb C$ and $q = e^{\frac {2\pi\sqrt{-1}}p}$, then $R$ is $q$-divisible (and therefore $q$-flat).
\item Assume that $q = 1$.
Then, quantum-flat means either that $R$ has no $\mathbb Z$-torsion (in which case $\mathrm{Char}(R) = 0$) or else that $R$ is an $\mathbb F_{p}$-algebra for some prime $p$ (and then $\mathrm{Char}(R) = p > 0$).
And quantum-divisible means that $R$ is an algebra over a field (whose characteristic is is the characteristic of $R$).
\item If $R = K[t]$ is a polynomial ring over a commutative ring $K$, and $q = t$, then $R$ is $q$-flat.
But $R$ is clearly \emph{not} $q$-divisible.
\item If $q = -1$, then $R$ is $q$-divisible because $(m)_{q}$ only takes values $0$ and $1$ (and same for $q = 0$).
\end{enumerate}
\end{xmp}

The flatness condition might sound odd but the quantum characteristic can have a rather strange behavior in general as the following examples show:

\begin{xmp}
\begin{enumerate}
\item If $q$ is the image of $X$ in $R = \mathbb Q[X]/(X^2 -1)$ then $q$ is a primitive square root of unity but $q\mathrm{-char}(R) = \mathrm{Char}(R) = 0$.
\item If $q$ is the image of $X$ in $R = \mathbf F_{2}[X]/(X^2 -1)$ then $q$ is a primitive square root of unity but $q\mathrm{-char}(R) = 4$ and $\mathrm{Char}(R) = 2$.
\item If $q= 1$ and $R = \prod \mathbb Z/n\mathbb Z$ then $q\mathrm{-char}(R) = \mathrm{Char}(R) = 0$ but $R$ has $(m)_{q}$-torsion for all $m\neq 0$.
\end{enumerate}
\end{xmp}

However, things get better under quantum flatness hypothesis:

\begin{prop}
Assume that $R$ is $q$-flat and let $p$ be a positive integer.
Then $q\mathrm{-char}(R) = p$ if and only if one of the following conditions is satisfied.
\begin{enumerate}
\item $q$ is a non-trivial primitive $p$-th root of unity.
\item $q = 1$ and $\mathrm{Char} (R) = p$.
\end{enumerate}
In both cases, $q$ is a $p$-th root of unity and, in particular, it is invertible.
\end{prop}

\begin{pf}
It follows from lemma \ref{pqzer} that $(p)_{q} = 0$ and that either $q$ is a non-trivial root of unity or else that $q = 1$ and $\mathrm{Char} (R) > 0$.
In the first case, $q$ is a primitive $p$-th root of unity if and only if $p$ is the smaller positive integer $m$ such that $q$ is an $m$-th root of unity.
In the second case, $\mathrm{Char} (R) = p$ if and only if $p$ is the smaller positive integer $m$ such that $\mathrm{Char} (R) \mid m$.
Therefore, the assertion follows from lemma \ref{pqzer} and the very definition of quantum characteristic.
$\quad \Box$
\end{pf}

When the quantum characteristic is even, we can do a little better:

\begin{prop} \label{even}
If $R$ is $q$-flat and $q\mathrm{-char}(R) = 2k > 0$, then $q^k = -1$.
\end{prop}

\begin{pf}
Since
$
(k)_{q}(2)_{q^k} = (2k)_{q} = 0
$
 and $R$ is $q$-flat, we must have $1 + q^k = (2)_{q^k} = 0$.
 $\quad \Box$
 \end{pf}
 
\begin{rmk}
The fact that $R$ is $q$-flat is crucial as the following example shows: if $R = \mathbb Z/8\mathbb Z$ and $q = 3$, we have $q\mathrm{-char}(R) = 4$ but $q^2 = 1 \neq -1$.
\end{rmk}

\begin{prop}
If $R$ has no $\mathbb Z$-torsion and $q\mathrm{-char}(R) = p > 0$, then $q$ is a primitive $p$-th root of unity.
\end{prop}

\begin{pf}
Since $(p)_{q} = 0$, lemma \ref{pqzer} tells us that $q$ is a $p$-th root of unity.
Assume that $q$ is not primitive. Then, there exists $1 \leq m < p$ with $p = mn$ such that $q^m = 1$.
It follows from \eqref{formul2} that
\begin{displaymath}
0 = (p)_{q} = (m)_{q}(n)_{1} = n (m)_{q}.
\end{displaymath}
Since $R$ has no $\mathbb Z$-torsion, necessarily $(m)_{q} = 0$ and this contradicts the minimality of $p$.
$\quad \Box$
\end{pf}

\begin{prop} \label{prim}
If $q\mathrm{-char}(R)$ is a prime number $p$, then $R$ is $q$-divisible (and therefore also $q$-flat).
\end{prop}

\begin{pf}
It follows from proposition \ref{divp} that $(m)_{q} = 0$ when $m$ is a multiple of $p$ and that $(m)_{q}$ is invertible otherwise.
$\quad \Box$
\end{pf}

\begin{rmk}
The condition $m$ prime to $p$ in the second statement of proposition \ref{divp} is necessary even if $R$ is $q$-flat as the following example shows.
If $R = \mathbb Z[\sqrt{-1}]$ and $q = \sqrt{-1}$, we have $q\mathrm{-char}(R)= 4$ and $(2)_{q} = 1 + \sqrt{-1} \not \in R^\times$.
\end{rmk}

\begin{lem} \label{moeb}
If $\chi_{m} \in \mathbb Z[t]$ denotes the $m$-th cyclotomic polynomial, we have whenever $n > 0$,
\begin{displaymath}
(n)_{q} = \prod_{m  \mid  n,m\neq 1} \chi_{m}(q).
\end{displaymath}
\end{lem}

\begin{pf}
We may assume that $R = \mathbb Z[t]$ and that $q = t$.
Then our assertion follows from the classical formula
\begin{displaymath}
1 - t^n = \prod_{n=md} \chi_{m}. \quad \Box
\end{displaymath}
\end{pf}

When $R$ is $q$-flat, the next result may be used to reduce some proofs to the case $R = \mathbb C$ and $q = e^{\frac{2\pi\sqrt{-1}}p}$:

\begin{prop} \label{rootp}
Assume that $R$ is $q$-flat with $q\mathrm{-char}(R) = p > 0$ and let $\zeta \in \mathbb Q^{\mathrm{alg}}$ be a primitive $p$-th root of unity.
Then, there exists a unique ring homomorphism $\mathbb Z[\zeta] \to R$ that sends $\zeta$ to $q$.
\end{prop}

\begin{pf}
Let us consider the unique ring homomorphism $u : \mathbb Z[t] \to R$ that sends $t$ to $q$.
With the notations of lemma \ref{moeb} we see that if $1 < n < p$, then $R$ has no $\chi_{n}(q)$-torsion (use the formula).
The same formula applied to the case $n = p$ then implies that $\chi_{p}(q) = 0$.
It follows that $\ker u$ contains the cyclotomic polynomial $\chi_{p}$ and factors therefore through $\mathbb Z[\zeta] := \mathbb Z[t]/\chi_{p}$. $\quad \Box$
\end{pf}

It will be quite important to understand the behavior of quantum characteristic, quantum flatness and quantum divisibility under the rising of $q$ to some power.

\begin{prop}
Assume $q\mathrm{-char}(R) = p > 0$ and let $k \in \mathbb N$ be such that $p \nmid k$ and $R$ has no $(k)_{q}$-torsion.
\begin{enumerate}
\item 
If $R$ is $q$-flat (resp. $q$-divisible), then $R$ is $q^k$-flat (resp. $q^k$-divisible).
\item
If $d = (p,k)$ denotes the greatest common divisor of $p$ and $k$, then $q^k\mathrm{-char}(R) = p/d$.
\end{enumerate}
\end{prop}

Note that the condition $p \nmid k$ is equivalent to $(k)_{q} \neq 0$ and that both hypothesis on $k$ are satisfied when $(k)_{q} \in R^\times$.

\begin{pf}
We let $m \in \mathbb N$.
Recall that $(km)_{q} = (k)_{q}(m)_{q^k}$.
Since we assume that $R$ has no $(k)_{q}$-torsion and that $(k)_{q} \neq 0$, we see that $(km)_{q} = 0$ is equivalent to $(m)_{q^k} = 0$.
But $(km)_{q} = 0$ means exactly that $p \mid km$ and this happens if and only if $p/d \mid m$.
Thus, we obtain the expected formula for the $q^k$-characteristic.

Now, we let $a \in R$ with $a \neq 0$.
If $R$ is $q$-flat and $(m)_{q^k}a = 0$, then we will have $(km)_{q}a = 0$ which implies that $(km)_{q} = 0$.
And we just saw that $(km)_{q} = 0$ if and only if $(m)_{q^k} = 0$.
Thus we see that $R$ is $q^k$-flat.

Assume now that $R$ is $q$-divisible.
We know that $(m)_{q^k} \neq  0$ if and only if $(km)_{q} \neq 0$, but then necessarily $(km)_{q} \in R^\times$ and therefore also $(m)_{q^k} \in R^\times$ because of the above equality $(km)_{q} = (k)_{q}(m)_{q^k}$ again.
And we see that $R$ is $q^k$-divisible.
$\quad \Box$
\end{pf}

\begin{rmk}
\begin{enumerate}
\item
The condition $p \nmid k$ in the proposition is really necessary because otherwise $R$ might be $q$-divisible but not even $q^k$-flat.
This is the case for example if $R = \mathbb Z/4 \mathbb Z$, $q = -1$ and $k = 2$.
\item Note also that if $p \mid k$, then $q^k\mathrm{-char} (R)$ is the usual characteristic of $R$.
In particular, it may be equal to $0$ whatever $p$ is.
\item Finally, the converse implications are false in general: if $R = \mathbb Z[\sqrt{-1}]$, $q = \sqrt{-1}$ and $k = 2$, we see that $R$ is $q^k$-divisible but not $q$-divisible.
\end{enumerate}
\end{rmk}

The last result of this section shows the relation between the dynamics of affine endomorphisms and quantum numbers.

\begin{prop} \label{sigmaen}
Assume $R$ is a commutative ring.
Let $A$ be an $R$-algebra and $\sigma$ an $R$-endomorphism of $A$.
Assume that $\sigma(x) = qx + h$ with $q,h \in R$.
Then, for all $n \in \mathbb N$ (or even $n \in \mathbb Z$ when $q \in R^\times$ and $\sigma$ is bijective), we have
\begin{enumerate}
\item $\sigma^n(x) = q^nx + (n)_{q}h$ 
\item $x - \sigma^n(x) = (n)_{q}(x - \sigma(x))$
\end{enumerate}
\end{prop}

\begin{pf}
By induction on $n \in \mathbb N$, we have
\begin{displaymath}
\sigma^{n+1}(x) = \sigma(q^nx + (n)_{q}h) = q^n(qx +h) + (n)_{q}h)
\end{displaymath}
\begin{displaymath}
= q^{n+1}x + ((n)_{q} + q^n)h = q^{n+1}x + (n+1)_{q}h
\end{displaymath}
and the case of a non-negative integer is settled.
Moreover, it follows that, when $q \in R^\times$ and $\sigma$ is bijective, we have
\begin{displaymath}
x = \sigma^{-n}(q^nx + (n)_{q}h) = q^n\sigma^{-n}(x) + (n)_{q}h
\end{displaymath}
and therefore,
\begin{displaymath}
\sigma^{-n}(x) = q^{-n}x - q^{-n}(n)_{q}h = q^{-n}x + (-n)_{q}h.
\end{displaymath}
It remains to prove the second assertion.
We have:
\begin{displaymath}
x - \sigma^n(x) = x - q^nx - (n)_{q}h  = (1 - q^n)x - (n)_{q}h = (n)_{q}(1-q)x - (n)_{q}h
\end{displaymath}
\begin{displaymath}
= (n)_{q}(x- qx - h) = (n)_{q}(x -\sigma(x)).
\quad \Box
\end{displaymath}
\end{pf}

\begin{rmk}
Even if we are only interested in commutative rings, non commutative ones might show up.
This is the case for example if $R$ is commutative, $A = R[x]$ denotes the polynomial ring in the variable $x$ over $R$ and we consider the ring $S$ of $R$-endomorphisms of $A$.
In particular, if $q \in R$, there exists a unique $\sigma \in S$ such that $\sigma(x) = qx$.
Then, in this case, we have $\sigma\mathrm{-char}(S) = q\mathrm{-char}(R)$.
\end{rmk}

\section{Quantum binomial coefficients} \label{binomial}

Recall that we work over a fixed ring $R$ and with a fixed $q \in R$.

\begin{dfn}
The \emph{$q$-factorial} (or \emph{quantum factorial}) of $m \in \mathbb N$ is
\begin{displaymath}
(m)_{q}! := \prod_{i=0}^{m-1}(m-i)_{q}.
\end{displaymath}
\end{dfn}

Of course, since $(1)_{q} = 1$, we could stop at $i = m-2$ as well.

In other words, we have
\begin{displaymath}
(0)_{q}! = (1)_{q}! = 1, \quad (2)_{q}! = (2)_{q} = 1 + q, \quad (3)_{q}! = (3)_{q}(2)_{q} = 1 + 2q + 2q^2 + q^3,
\end{displaymath}
and for bigger $m$,
\begin{displaymath}
(m)_{q}! := (m)_{q} (m-1)_{q} \cdots (3)_{q}(2)_{q}.
\end{displaymath}

\begin{xmp}
\begin{enumerate}
\item If $q = 1$ and $\mathbb Z \subset R$, then $(m)_{q}! = m!$ is the usual factorial.
\item More generally, when $q = 1$, we have $(m)_{q}! = m!1_{R}$.
In particular, we see that $(m)_{q}! = 0$ for $m \geq p$ when $q = 1$ and $\mathrm{Char}(R) = p > 0$.
\item When $R = \mathbb Q(t)$ and $q = t$, we have
\begin{displaymath}
(m)_{q}! = \frac {(1 -t^m)(1-t^{m-1})\cdots (1 - t^2)(1-t)}{(1 -t)^m}
\end{displaymath}
\item If $R = \mathbb C$ and $q = e^{2\pi\sqrt{-1}/p}$ with $p$ an integer $\geq 2$, we have
\begin{displaymath}
(m)_{q}! = \left\{
\begin{array}{ccc}
\frac {(1 -  e^{2m\pi\sqrt{-1}/p})(1- e^{2(m-1)\pi\sqrt{-1}/p})\cdots (1 -  e^{2\pi\sqrt{-1}/p})}{(1 - e^{2\pi\sqrt{-1}/p})^m} & \mathrm{if} & m < p
\\ 0 & \mathrm{if} & m \geq p.
\end{array}\right.
\end{displaymath}
\end{enumerate}
\end{xmp}

\begin{prop}
If $q\mathrm{-char}(R) = p$, then $(m)_{q}! = 0$ for $m \geq p$.
\end{prop}

\begin{pf} Immediate consequence of the definition. $\quad \Box$
\end{pf}

\begin{prop}
For all $m \in \mathbb N$, we have
\begin{displaymath}
(1 -q)^m(m)_{q}! = \prod_{i=0}^{m-1} (1 - q^{m-i}).
\end{displaymath}
In particular, if $1 - q$ is invertible in $R$, we have
\begin{displaymath}
(m)_{q}! = \frac {(1 -q^m)(1-q^{m-1})\cdots (1 - q^2)(1-q)}{(1 -q)^m}.
\end{displaymath}
\end{prop}

\begin{pf}
Follows from lemma \ref{explicit}.$\quad \Box$
\end{pf}

\begin{prop} \label{moeb2}
If $\chi_{m} \in \mathbb Z[t]$ denotes the $m$-th cyclotomic polynomial, we have for all $n \in \mathbb N$,
\begin{displaymath}
(n)_{q}! = \prod_{m\neq1} \chi_{m}(q)^{[\frac nm]}
\end{displaymath}
where $[\frac nm]$ is the integer part of $\frac nm$.
\end{prop}

\begin{pf}
We saw in lemma \ref{moeb} that we have for all $k \in \mathbb N$,
\begin{displaymath}
(k)_{q} = \prod_{m  \mid  k,m\neq 1} \chi_{m}(q).
\end{displaymath}
On the other hand, one easily sees that for all $m \in \mathbb N$, we have
\begin{displaymath}
\# \{k \leq n, m \mid k\} = \left[\frac nm\right],
\end{displaymath}
and the formula follows.
 $\quad \Box$
\end{pf}

\begin{dfn}
The \emph{$q$-binomial coefficients} (or \emph{quantum binomial coefficients}) are defined by induction for $k,n \in \mathbb N$ via \emph{Pascal identities}
\begin{displaymath}
{n \choose k}_{q} = {n-1 \choose k-1}_{q} + q^k{n-1 \choose k}_{q}
\end{displaymath}
with
\begin{displaymath}
{0 \choose k}_{q} = \left\{\begin{array}{l} 1\ \mathrm{if}\ k = 0 \\ 0 \ \mathrm{otherwise}. \end{array}\right. \end{displaymath}
\end{dfn}

\begin{rmk}
If we want to use the ``symmetric quantum state'' (as it is usually the case in quantum group theory),
\begin{displaymath}
[n]_{v} := \frac {v^n -v^{-n}}{v - v^{-1}},
\end{displaymath}
then we will have
\begin{displaymath}
[n]_{v} ! = \frac {1}{v^{\frac {n(n-1)}2}} (n)_{v^2}! \quad \mathrm{and} \quad \left[{n \atop k}\right]_{v} = \frac {1}{v^{\frac {n(n-k)}2}}{n \choose k}_{v^2}.
\end{displaymath}
\end{rmk}

\begin{xmp}
When $R = \mathbb Z$ and $q$ is a power of a prime $p$, then ${n \choose k}_{q}$ is the number of rational points of the Grassmanian $\mathbb{G}(n,k,q)$.
Said differently, this is the number of vector subspaces of dimension $k$ in a vector space of dimension $n$ over a field with $q$ elements.
This is easily checked (see also \cite{KacCheung02}, Theorem 7.1).
\end{xmp}

\begin{prop}
We  have for all $n,k \in \mathbb N$,
\begin{displaymath}
{n \choose k}_{q}  \prod_{i=0}^{k-1}(k-i)_{q} = \prod_{i=0}^{k-1}(n-i)_{q}.
\end{displaymath}
In particular, if $q\mathrm{-char}(R) = 0$ and $R$ is $q$-divisible, then
\begin{equation} \label{bino}
{n \choose k}_{q} = \frac {(n)_{q}!}{(k)_{q}!(n-k)_{q}! }.
\end{equation}
\end{prop}

\begin{pf}
In order to prove the first assertion, we may first assume that $R = \mathbb Z[t]$ and $q = t$, and then specialize to any $R$ and $q$.
We may even assume that $R = \mathbb Q(t)$.
In particular, all non zero $q$-integers will be invertible in $R$ and it is therefore sufficient to prove the second assertion.
We can use lemma \ref{explicit} in order to show that the right member of equality \eqref{bino} satisfies the induction property of the left member.
This works as follows:
\begin{displaymath}
\frac {(n-1)_{q}!}{(k-1)_{q}!(n-k)_{q}! } + q^k\frac {(n-1)_{q}!}{(k)_{q}!(n-k-1)_{q}! }
\end{displaymath}
\begin{displaymath}
= \frac {((k)_q + q^k(n-k)_q)(n-1)_{q}!} {(k)_{q}!(n-k)_{q}! } = \frac {(n)_{q}!}{(k)_{q}!(n-k)_{q}! }. \quad \Box
\end{displaymath}
\end{pf}

\begin{cor}
We  have for all $n,k \in \mathbb N$,
\begin{displaymath}
{n \choose n-k}_{q}  = {n \choose k}_{q}.
\end{displaymath}
\end{cor}

\begin{pf}
We may assume as above than $R = \mathbb Q(t)$ and $q = t$ and use formula \eqref{bino}.
$\quad \Box$
\end{pf}

\begin{cor}
We have for all $k,n \in \mathbb N$,
\begin{displaymath}
 (1 - q^k)(1 - q^{k-1}) \cdots (1 - q){n \choose k}_{q}  = (1 - q^n)(1 - q^{n-1}) \cdots (1 -q^{n-k+1}).
\end{displaymath}

In particular, if $1 - q^i$ is invertible for all $0 < i \leq k$, we will have
\begin{displaymath}
{n \choose k}_{q}  = \frac {(1 - q^n)(1 - q^{n-1}) \cdots (1 -q^{n-k+1})}{(1 - q^k)(1 - q^{k-1}) \cdots (1 - q)}. \quad \Box
\end{displaymath}
\end{cor}

\begin{rmk}
\begin{enumerate}
\item The rational function
\begin{displaymath}
\frac {(1 - t^n)(1 - t^{n-1}) \cdots (1 -t^{n-k+1})}{(1 - t^k)(1 - t^{k-1}) \cdots (1 - t)} \in \mathbb Q(t)
\end{displaymath}
actually lives in $\mathbb Z[t]$ and is called a \emph{Gaussian polynomial}.
It is identical to the binomial coefficient ${n \choose k}_{t}$.
\item 
One may prove many properties of quantum binomial coefficients by reducing to the case $R = \mathbb Q(t)$ and $q = t$ and using various references in the literature (see for example \cite{Kassel95}, section IV.2).
\item Actually, one may as well assume that $R = \mathbb C$ because it is always possible to embed $\mathbb Q(t)$ into $\mathbb C$ by sending $t$ to any transcendental $q \in \mathbb C$.
\end{enumerate}
\end{rmk}

One may also define the quantum binomial coefficients as a product as we can see for example in \cite{KnuthWilf89}:

\begin{cor} \label{moeb3}
We have for all $n \in \mathbb N$,
\begin{displaymath}
{n \choose k}_{q} = \prod_{[\frac nm]> [\frac km]+[\frac {n-k}m]} \chi_{m}(q)
\end{displaymath}
where $\chi_{m} \in \mathbb Z[t]$ denotes the $m$-th cyclotomic polynomial and $[r]$ denotes the integer part of a real number $r$.
\end{cor}

The condition under the product says that the sum of the rests in the euclidean division of $k$ and $n-k$ by $m$ is at least equal to $m$.

\begin{pf}
We may assume that $q\mathrm{-char}(R) = 0$ and $R$ is $q$-divisible.
Then, formula \eqref{bino} and proposition \ref{moeb2} give
\begin{displaymath}
{n \choose k}_{q}  = \prod_{m\neq 1} \chi_{m}(q)^{[\frac nm]-[\frac km]-[\frac {n-k}m]}
\end{displaymath}
and we have $[\frac nm]-[\frac km]-[\frac {n-k}m] = 0$ unless $[\frac nm]> [\frac km]+ [\frac {n-k}m]$ in which case $[\frac nm]-[\frac km]-[\frac {n-k}m] = 1$.
Note that this never happens when $m = 1$.
$\quad \Box$
\end{pf}

\begin{prop} \label{transi}
We have
\begin{displaymath}
\forall n,j,k \in \mathbb N, \quad {n \choose j}_q{j \choose k}_q = {n \choose k}_q{n -k  \choose n -j}_q.
\end{displaymath}
\end{prop}

\begin{pf}
We may assume that $R = \mathbb Q(t)$ and $q = t$ and our formula reads
\begin{displaymath}
 \frac {(n)_{t}!}{(j)_{t}!(n-j)_{t}! } \quad \frac {(j)_{t}!}{(k)_{t}!(j-k)_{t}! } =  \frac {(n)_{t}!}{(k)_{t}!(n-k)_{t}! } \quad  \frac {(n-k)_{t}!}{(n-j)_{t}!(j-k)_{t}! }. \quad \Box
\end{displaymath}
\end{pf}

We can also state and prove the \emph{quantum Chu-Vandermonde identity}:

\begin{lem} \label{binobin}
We have
\begin{displaymath}
\forall n,m,k \in \mathbb N, \quad  {n+m \choose k}_{q} = \sum_{i=0}^k  q^{i(m-k+i)} {n \choose i}_{q} {m \choose k-i}_{q}.
\end{displaymath}
\end{lem}

Recall that, with our conventions, we have ${n \choose i}_{q} = 0$ for $i > n$ and ${m \choose k-i}_{q} = 0$ for $k-i > m$.

\begin{pf} This is shown to be true by induction on $m$.
We will have
\begin{displaymath}
 {n+m \choose k}_{q} =  {n+m-1 \choose k-1}_{q} +  q^k{n+m-1 \choose k}_{q} 
\end{displaymath}
\begin{displaymath}
= \sum_{i=0}^{k-1}  q^{i(m-k+i)}{n \choose i}_{q} {m-1 \choose k-1-i}_{q} + q^k \sum_{i=0}^k  q^{i(m-1-k+i)}{n \choose i}_{q} {m-1 \choose k-i}_{q}
\end{displaymath}
\begin{displaymath}
= \sum_{i=0}^k  q^{i(m-k+i)}{n \choose i}_{q} \left({m-1 \choose k-1-i}_{q} + q^{k-i}  {m-1 \choose k-i}_{q} \right) 
\end{displaymath}
\begin{displaymath}
= \sum_{i=0}^k  q^{i(m-k+i)}{n \choose i}_{q} {m \choose k-i}_{q}. \quad \Box
\end{displaymath}
\end{pf}

\begin{lem} \label{qbin}
Assume $q\mathrm{-char}(R) = p > 0$ and $R$ is $q$-flat, then
\begin{displaymath}
 {p \choose k}_{q} =
\left\{\begin{array}{clc}   1 & \mathrm{if}\ k = 0\ \mathrm{or}\ k = p \\ 0 & \mathrm{otherwise.} & \quad \Box \end{array}\right. 
\end{displaymath}
\end{lem}

\begin{pf}
We may assume $0 < k < p$.
Since $q\mathrm{-char}(R) = p > 0$, we will have
\begin{displaymath}
 (k)_{q}(k-1)_{q} \cdots (2)_{q}{p \choose k}_{q} = (p)_{q}(p-1)_{q} \cdots (p-k+1)_{q} = 0.
\end{displaymath}
And since we assume that $R$ is $q$-flat, we must have
\begin{displaymath}
{p \choose k}_{q} = 0. \quad \Box
\end{displaymath}
\end{pf}

\begin{rmk}
The condition will always be satisfied when $p$ is prime.
Actually, when $q = 1$, the flatness condition is equivalent to $p$ being prime.
However, this is not necessary in general.
\end{rmk}

\begin{xmp}
\begin{enumerate}
\item Assume $R = \mathbb C$ and $q = e^{\frac {2\sqrt{-1}\pi}{p}}$ with $p \in \mathbb N$ (not necessary prime) and $p \geq 2$.
Then we have ${p \choose k}_{q} = 0$ for $0 < k < p$.
\item Assume $R$ is an $\mathbb F_{p}$-algebra for some prime number $p$.
Then ${p \choose k} = 0$ for $0 < k < p$.
\item Assume $R = \mathbb Z/4 \mathbb Z$ and $q = 1$. Then we have ${4 \choose 2}_{q} = 2 \neq 0$.
\end{enumerate}
\end{xmp}

We can now prove the \emph{quantum Lucas theorem} (see also lemma 24.1.2 of \cite{Lusztig10}):

\begin{prop} \label{binexp}
Assume $q\mathrm{-char}(R) = p > 0$ and $R$ is $q$-flat.
Let $n,k,i,j \in \mathbb N$ with $i,j < p$.
Then, we have
\begin{displaymath}
 {np + i \choose kp+j}_{q} = {n \choose k}{ i \choose j}_{q}.
\end{displaymath}
With our convention, it means in particular that
\begin{displaymath}
{np+i \choose kp+j}_{q} = 0 \quad \mathrm{if}\ 0 \leq i < j < p.
\end{displaymath}
\end{prop}

\begin{pf}
We proceed by induction on $n$ and $i$, and we use the quantum Pascal identity
\begin{equation} \label{indfor}
{np+i \choose kp+j}_{q} = {np+i-1\choose kp+j-1}_{q}  + q^{kp+j} {np+i-1\choose kp+j}_{q}.
\end{equation}
We only do the non trivial cases.

Assume first that $n,k >0$ but $i = j = 0$.
Then the formula reads
\begin{displaymath}
{np \choose kp}_{q} = {(n-1)p+p-1 \choose (k-1)p + p-1}_{q}  + q^{kp} {(n-1)p+p-1\choose kp}_{q}
\end{displaymath}
\begin{displaymath}
= {n-1 \choose k-1}{p-1 \choose p-1}_{q}  +  {n-1 \choose k} {p-1\choose 0}_{q} =  {n-1 \choose k-1}  +  {n-1 \choose k}  = {n \choose k}
\end{displaymath}
as expected.

Assume now that $n,j > 0$ but $i = 0$.
Then formula \eqref{indfor} reads
\begin{displaymath}
{np \choose kp+j}_{q} = {(n-1)p+p-1 \choose kp+j-1}_{q}  + q^{kp+j} {(n-1)p+p-1\choose kp+j}_{q}
\end{displaymath}
\begin{displaymath}
= {n-1 \choose k}{p-1 \choose j-1}_{q}  + q^{j} {n-1 \choose k} {p-1\choose j}_{q} = {n-1 \choose k} \left({p-1 \choose j-1}_{q}  + q^{j} {p-1\choose j}_{q} \right)
\end{displaymath}
\begin{displaymath}
= {n-1 \choose k} {p \choose j}_{q} = 0
\end{displaymath}
thanks to lemma \ref{qbin}.

Now, if $i,j > 0$, the formula reads
\begin{displaymath}
{np+i \choose kp+j}_{q} = {n \choose k}{ i-1 \choose j-1}_{q} + q^j {n \choose k}{ i-1 \choose j}_{q}
\end{displaymath}
\begin{displaymath}
= {n \choose k}\left({ i-1 \choose j-1}_{q} + q^j { i-1 \choose j}_{q}\right) = {n \choose k}{ i \choose j}_{q}.
\end{displaymath}

Finally, in the case $i,k >0$ but $j = 0$, formula \eqref{indfor} gives
\begin{displaymath}
{np+i \choose kp}_{q} = {np+i-1\choose kp-1}_{q}  + q^{kp} {np+i-1 \choose kp}_{q}
\end{displaymath}
\begin{displaymath}
= {np+i-1\choose (k-1)p + p-1}_{q}  +  {np+i-1 \choose kp}_{q} = {n \choose k-1}{ i-1 \choose p-1}_{q} + {n \choose k}{ i-1 \choose 0}_{q}
\end{displaymath}
\begin{displaymath}
= {n \choose k-1} \times 0 + {n \choose k} \times 1 = {n \choose k} = {n \choose k}{ i \choose 0}_{q}
\end{displaymath}
because $0 \leq i- 1 < p-1 < p$.
$\quad \Box$
\end{pf}

\begin{rmk}
\begin{enumerate}
\item We recover the usual Lucas theorem in arithmetics from the case $R = \mathbb F_{p}$ and $q = 1$ of the proposition:
if $n = \sum a_{i}p^i$ and $k = \sum b_{i}p^{i}$ denote the $p$-adic expansions of $n$ and $k$ ($p$ a \emph{prime} number), we have
\begin{displaymath}
{n \choose k} \equiv \prod_{i} {a_{i} \choose b_{i}} \mod p.
\end{displaymath}
\item In the case where $R = \mathbb C$ is the field of complex numbers and $q = \zeta$ is a primitive $p$-th root of unity, we recover proposition 2.1 of \cite{GuoZeng06}.
\item Using proposition \ref{rootp}, one can also derive the quantum Lucas theorem from the case $R = \mathbb C$.
This gives a proof of the classical Lucas theorem using the theory of complex functions.
\end{enumerate}
\end{rmk}

Finally, we prove the binomial quantum formula:

\begin{prop} \label{bincof}
Assume that $R$ is commutative and let $A$ be a commutative $R$-algebra.
Then, we have
\begin{displaymath}
\forall n \in \mathbb N, \quad \prod_{i=0}^{n-1}(q^ix+y) = \sum_{k=0}^n q^{\frac {k(k-1)}2}{n \choose k}_q x^{k}y^{n-k}.
\end{displaymath}
\end{prop}

\begin{pf}
By induction on $n$, we see that
\begin{displaymath}
\prod_{i=0}^{n-1}(q^ix+y) = \left(\sum_{k=0}^{n-1} q^{\frac {k(k-1)}2}{n-1 \choose k}_q x^{k}y^{n-1-k} \right)(q^{n-1}x + y)
\end{displaymath}
\begin{displaymath}
= \left(\sum_{k=0}^{n-1} q^{\frac {k(k-1)}2+n-1}{n-1 \choose k}_q x^{k+1}y^{n-1-k} \right)
+ \left(\sum_{k=0}^{n-1} q^{\frac {k(k-1)}2}{n-1 \choose k}_q x^{k}y^{n-k} \right)
\end{displaymath}
\begin{displaymath}
= \left(\sum_{k=1}^{n} q^{\frac {(k-1)(k-2)}2+n-1}{n-1 \choose k-1}_q x^{k}y^{n-k} \right) + \left(\sum_{k=0}^{n-1} q^{\frac {k(k-1)}2}{n-1 \choose k}_q x^{k}y^{n-k} \right)
\end{displaymath}
\begin{displaymath}
= \left(\sum_{k=0}^{n}q^{\frac {k(k-1)}2} \left(q^{n-k}{n-1 \choose k-1}_q + {n-1 \choose k}_q\right) x^{k}y^{n-k}  \right)
\end{displaymath}
\begin{displaymath}
= \sum_{k=0}^n q^{\frac {k(k-1)}2}{n \choose k}_q x^{k}y^{n-k} . \quad \Box
\end{displaymath}
\end{pf}

\begin{rmk} 
If instead of assuming $A$ commutative, we make the supposition that $yx = qxy$ (\emph{quantum plane identity}), then the formula becomes even nicer:
\begin{displaymath}
\forall n \in \mathbb N, \quad (x+y)^n = \sum_{k=0}^n {n \choose k}_q x^{k}y^{n-k}
\end{displaymath}
(see proposition IV.2.2 of \cite{Kassel95} for example).
\end{rmk}

\section{Quantum rational numbers} \label{rational}

We are going to define the quantum states of a rational number.
We might call them quantum rational numbers (as in \cite{Nathanson06}) but they should not be confused with the quantum rational numbers that appear in quantum physics (see in chapter 6 of \cite{Saller06} for example).

We start with some generalities about roots in monoids, generalizing divisibility in additive (commutative) monoids.
In the end, we will apply these considerations to the multiplicative monoid of $R$.

We recall that a monoid $S$ is a set endowed with a law which is associative with unit.
Usually, this law is written multiplicatively, but we might also use the addition when the law is commutative.

\begin{dfn} \label{monroot}
Let $S$ be a  monoid.
A family $\{s_n\}_{n \in D}$ with $\emptyset \neq D \subset \mathbb N \setminus \{0\}$ is a \emph{system of roots} in $S$ if it satisfies:
\begin{displaymath}
\forall n,n' \in D, \forall m,m' \in \mathbb N, \quad m'n=mn'  \Rightarrow s_{n'}^{m'} = s_n^{m}.
\end{displaymath}
\end{dfn}

In other words, we require that $s_{n}^m$ only depends on $r := \frac mn \in \mathbb Q$ when $n \in D$ and $m \in \mathbb N$.
In particular, $s := s_n^{n}$ does not depend on $n \in D$ and we will also call $\{s_n\}_{n \in D}$  a \emph{system of roots of $s$}.

We specialize a little bit the definition:

\begin{dfn} \label{monroot2}
Let $S$ be a monoid and $\underline s := \{s_n\}_{n \in D}$ a system of roots of $s \in S$.
\begin{enumerate}
\item In the case $D := \{p\}$, we will call $s_p$ a \emph{$p$-th root of $s$}.
\item In the case $D := \{p^{i}, i \in \mathbb N\}$, we will call $\underline s$  a \emph{system of $p$-th roots of $s$}.
\item In the case $D = \mathbb N \setminus \{0\}$, we will call $\underline s$ a \emph{complete system of roots of $s$}.
\end{enumerate}
\end{dfn}

For $p$-th roots, or more generally, for systems of $p$-th roots, there exists a simpler alternative definition:

\begin{prop}
Let $S$ be a monoid and $s \in S$.
\begin{enumerate}
\item If $p \in \mathbb N \setminus \{0\}$, giving a $p$-th root of $s$ is equivalent to giving an element $s_1 \in S$ such that $s = s_1^p$.
\item If $p \in \mathbb N \setminus \{0\}$, giving a system of $p$-th roots of $s$ is equivalent to giving a sequence $\{s_{i}\}_{i \in \mathbb N}$ of $s_i \in S$ such that $s_0 = s$ and $s_{i+1}^p = s_i$.
\end{enumerate}
\end{prop}

\begin{pf}
In the first assertion, the condition of the definition is void and we have $s_p^p=s$.
Thus, we obtain the result after a renumbering $s_p \leadsto s_1$.

For the second assertion, one easily checks that the condition in the definition is implied by
\begin{displaymath}
\forall i \in \mathbb N, \quad s_{p^{i+1}}^p = s_{p^i}.
\end{displaymath}
And the result therefore follows also from a renumbering $s_{p^{i}} \leadsto s_i$.
$\quad \Box$
\end{pf}

Any monoid $S$ has a \emph{natural preorder} (reflexive and transitive relation) given by
\begin{displaymath}
s \leq s' \quad  \Leftrightarrow \quad \exists m \in \mathbb N, s' = s^m.
\end{displaymath}
For example, the natural preorder on the additive monoid $\mathbb N$ is given by divisibility (and not the usual order on $\mathbb N$).
Note that any morphism of monoids preserves the preorder.
Finally, recall that a preordered set is \emph{inductive} (or \emph{directed}) if any couple has an upper bound:
\begin{displaymath}
\forall s,s' \in S, \quad \exists s'' \in S \quad s \leq s'' \ \mathrm{and}\ s' \leq s''.
\end{displaymath}

\begin{rmk}
\begin{enumerate}
\item
When the index set $D$ is inductive (for divisibility), the condition of definition \ref{monroot} is equivalent to
\begin{displaymath}
\forall n,n' \in D, \forall m \in \mathbb N, \quad n = mn' \Rightarrow s_{n'} = s_{n}^{m}.
\end{displaymath}
\item When $D$ is inductive, the family $\{s_n\}_{n \in D}$ is inductive for the \emph{reverse} preorder.
\item
In the special cases above, the index set is inductive (and therefore, the system of roots is inductive for the reverse preorder).
\end{enumerate}
\end{rmk}

\begin{dfn}
Let $N$ be a submonoid of the additive monoid $\mathbb Q_{\geq0}$.
A denominator for $N$ is an element $n \in \mathbb N \setminus \{0\}$ such that $\frac 1n \in N$.
A subset $D \subset \mathbb N \setminus \{0\}$ is a \emph{full set of denominators} for $N$ if $\frac 1D := \{\frac 1n, n \in D\}$ is a set of generators for $N$.
\end{dfn}

If $E \subset N$ is a set of generators for an additive (commutative) monoid $N$, we will write $N = \mathbb NE$.
Thus, we see that $D$ is a full set of denominators for $N$ if $N = \mathbb N\frac 1D$.

Recall also that if $N$ is an additive monoid, there exists a smallest abelian group $\pm N$ that contains $N$.
More precisely, the forgetful functor from abelian groups to commutative monoids has a left adjoint $N \mapsto \pm N$.
Note that when $N$ is a submonoid of $\mathbb Q_{\geq 0}$, we may assume $\pm N \subset \mathbb Q$, and then we have $N = \pm N \cap \mathbb Q_{\geq 0}$.

\begin{xmp}
\begin{enumerate}
\item For $D := \{p\}$ with $p \in \mathbb N \setminus \{0\}$, we have
\begin{displaymath}
\mathbb N\frac 1D = \mathbb N\frac 1p := \{\frac mp, \quad m \in \mathbb N\}
\end{displaymath}
and
\begin{displaymath}
\pm \mathbb N\frac 1D = \mathbb Z\frac 1p := \{\frac mp, \quad m \in \mathbb Z\}.
\end{displaymath}
\item If $D := \{p^{i}, i \in \mathbb N\}$ with $p \in \mathbb N \setminus \{0\}$, then
\begin{displaymath}
\mathbb N\frac 1D = \mathbb N\left[\frac 1p\right] := \{r \in \mathbb Q_{\geq 0}, \exists i \in \mathbb N, p^ir \in \mathbb N\}
\end{displaymath}
and
\begin{displaymath}
\pm \mathbb N\frac 1D = \mathbb Z\left[\frac 1p\right] := \{r \in \mathbb Q, \exists i \in \mathbb N, p^ir \in \mathbb Z\}.
\end{displaymath}
\item If $D = \mathbb N \setminus \{0\}$, then $\mathbb N\frac 1D  = \mathbb Q_{\geq0}$ and $\mathbf \pm \mathbb N\frac 1D  = \mathbb Q$.
\end{enumerate}
\end{xmp}

\begin{lem} \label{stupid}
If $m, n$ are two denominators for a submonoid $N$ of $\mathbb Q_{\geq 0}$, then their least common multiple $p$ is also a denominator for $N$.
\end{lem}

\begin{pf}
We are given two denominators $m,n$ of $N$.
Let us denote by $d$ their greatest common divisor and by $p$ their least common multiple.
We can write $d = um + vn$ with $u, v \in \mathbb Z$ and it follows that $\frac 1p = \frac un + \frac vm \in \pm N$ and therefore $\frac 1p \in N$.
$\quad \Box$
\end{pf}

\begin{prop}
A submonoid $N \subset \mathbb Q_{\geq 0}$ has a full set of denominators if and only if $N = \{0\}$ or $1 \in N$.
If this is the case, it has a full inductive set of denominators.
Actually, if $1 \in N$, the set
\begin{displaymath}
D := \{n \in \mathbb N \setminus \{0\}, \frac 1n \in N\}
\end{displaymath}
of all denominators of $N$ is a full inductive set of denominators for $N$.
\end{prop}

Of course, the condition $1 \in N$ is equivalent to $\mathbb N \subset N$.

\begin{pf}
The condition is necessary.
More precisely, there exists $n \in D$ and we have $1 = n \times \frac 1n \in N$.
In order to check that the condition is also sufficient, we only have to prove the last assertion.

Let us assume that $1 \in N$.
If $r \in N$, we can write $r = \frac mn$ with $m,n \in \mathbb N$ coprime and $n \neq 0$.
Thus, there exists $u, v \in \mathbb Z$ with $um+vn=1$ and it follows that $\frac 1n = ur+v \in \pm N$.
Therefore, we can write $r = m \times \frac 1n$ with $m \in \mathbb N$ and $\frac 1n \in \pm N \cap \mathbb Q_{\geq 0} = N$.
It means that all denominators make a full set of denominators.
We still have to show that this is an inductive set.
Actually, this follows from lemma \ref{stupid}.
$\quad \Box$
\end{pf}

\begin{prop} \label{fulend}
\begin{enumerate}
\item If a submonoid $N \subset \mathbb Q_{\geq 0}$ has a finite full set of denominators $D$, then $N = \mathbb N\frac 1p$ for some $p \in \mathbb N$.
\item If $D$ is a full inductive set of denominators for a submonoid $N$ of $\mathbb Q_{\geq 0}$, then
\begin{displaymath}
N = \cup_{n \in D} \mathbb N\frac 1n \simeq \varinjlim_{n \in D} \mathbb N\frac 1n.
\end{displaymath}
\end{enumerate}\end{prop}

Note that the second assertion means that any $r \in N$ may be written on the form $r = \frac mn$ with $m \in \mathbb N$ and $n \in D$.

\begin{pf}
We prove the first assertion.
Let $D$ be a finite set of positive integers and $p$ the least common multiple of all elements of $D$.
Clearly, we have $N \subset \mathbb N\frac 1p$ and it only remains to check that $\frac 1p \in N$.
By induction, this  will easily follow from the case $D = \{m,n\}$.
And we can use lemma \ref{stupid}.

In order to prove the second assertion, it is sufficient to check that $\cup_{n \in D} \mathbb N\frac 1n$ is a submonoid of $\mathbb Q_{\geq 0}$.
But, since $D$ is inductive, if $n, n' \in D$, there exists $n'' \in D$ such that $n=dn''$ and $n'=d'n''$ with $d,d' \in \mathbb N$.
Therefore, if $m, n \in \mathbb N$, we have
\begin{displaymath}
\frac mn + \frac {m'}{n'} = \frac {md + m'd'}{n''} \in \mathbb N\frac 1{n''}. \quad \Box
\end{displaymath}
\end{pf}

\begin{prop} \label{monomap}
Let $N$ be a submonoid of the additive monoid $\mathbb Q_{\geq 0}$ that contains $\mathbb N$, $S$ a (multiplicative) monoid and $s \in S$.
\begin{enumerate}
\item If
\begin{equation} \label{maproot}
\xymatrix@R0cm{
N \ar[r] & S
\\ r \ar@{|->}[r] & s^r
}
\end{equation}
is a morphism of monoids that sends $1$ to $s$ and $D$ is a full set of denominators for $N$, then the sequence $\{s^{\frac 1n}\}_{n \in D}$ is a system of roots of $s$.
\item Conversely, if $D$ is a full inductive set of denominators for $N$ and $\{s_n\}_{n \in D}$ is a system of roots for $s$, there exists a unique map \eqref{maproot} with $s^{\frac 1n} = s_n$ for all $n \in D$.
\end{enumerate}
Moreover, the map \eqref{maproot} extends uniquely to $\pm N$ if and only if $s$ is invertible in $S$.
\end{prop}

\begin{pf}
Under the hypothesis of the first assertion, one easily checks that the conditions for a system of roots are satisfied.
More precisely, we have for all $n \in D$, $(s^{\frac 1n})^n = s$.
Moreover, if $r = \frac mn$ with $n \in D$, then $(s^{\frac 1n})^m = s^r$ will only depend on $r$.

Using proposition \ref{fulend} and uniqueness, the second assertion will follow from the case $D = \{p\}$ which in turn follows from the fact that $\mathbb N\frac 1p$ is isomorphic to $\mathbb N$ as an abstract monoid.

Finally, if $s$ is invertible and $\frac 1s$ denotes its inverse, one can extend the map \eqref{maproot} to $\pm N$ by sending $-r$ to $(\frac 1s)^r$.
Of course, one must check that this defines a morphism of monoids.
This is left to the reader.
Conversely, if the maps extends to $\pm N$, the image of $-1$ must be an inverse for $s$.
$\quad \Box$
\end{pf}

\begin{cor}
If we are given a system of $p$-th roots ${s_{p^{i}}}$ of $s$ for \emph{all} prime $p$, this will extend uniquely to a complete system of roots of $s$. 
\end{cor}

\begin{pf}
Uniqueness follows from the fact that $\sum_{p} \mathbb N\left[\frac 1p\right] = \mathbb Q_{\geq 0}$ (i.e. $\mathbb Q_{\geq 0}$ is the smallest submonoid containing all $\mathbb N\left[\frac 1p\right]$ for $p$ prime).
Existence follows from the fact that $\mathbb N[\frac 1{p_{1}}] \cap \mathbb N[\frac 1{p_{2}}] = \mathbb N$ for $p_{1}, p_{2}$ distinct primes.
Details are left to the reader.
$\quad \Box$
\end{pf}

As we said above, we want to apply the theory to the multiplicative monoid of our ring $R$ and the element $q$.

\begin{xmp}
For $R = \mathbb C$ and $q = \rho e^{\sqrt{-1}\theta}$, we can consider the morphism of groups
\begin{displaymath}
\xymatrix@R0cm{
\mathbb Q \ar[r] & \mathbb C^\times
\\ r \ar@{|->}[r] & \rho^r e^{\sqrt{-1}r\theta}.
}
\end{displaymath}
It provides us with a complete system of roots of $q$.
More generally, if $K$ is algebraically closed, there always exists a complete system of roots of $q \in K$.
\end{xmp}

Recall that if $K$ is a commutative ring, the forgetful functor from $K$-algebras to monoids has a left inverse.
Actually, if $N$ is an additive monoid, the associated $K$-algebra is the free module on the abstract basis $\{t^r\}_{r \in N}$ and multiplication is given by $t^{r_1}t^{r_2} = t^{r_1 + r_2}$.

When $N$ is a submonoid of $\mathbb Q_{\geq 0}$ with set of denominators $D$, we will denote the $K$-algebra of $N$ by $K[t^{\frac 1D}]$ (even if it actually only depends only on $N$ and not on $D$).
This is the ring of \emph{Puiseux polynomials} with denominators in $D$.
Note that
\begin{displaymath}
K[t^{\frac 1D}] = \varinjlim_{n \in D} K[t^{\frac 1n}]
\end{displaymath}
when $D$ is inductive.
Actually, the map
\begin{displaymath}
\xymatrix@R0cm{ K[t] \ar[r] & K[t^{\frac 1n}]
\\ t \ar@{|->}[r] & t^{\frac 1n}}
\end{displaymath}
is obviously an isomorphism and we could as well write
\begin{displaymath}
K[t^{\frac 1D}] = \varinjlim_{t \mapsto t^n,n \in D} K[t].
\end{displaymath}
We will also denote the $K$-algebra of $\pm N$ by $K[t^{\pm\frac 1D}]$ and, when $K$ is a field, we will denote by $K(T^{\frac 1D})$ the fraction field of $K[t^{\frac 1D}]$.

The $K$-algebra $K[t^{\frac 1D}]$ has the following universal property:

\begin{prop}
Assume that $R$ is a $K$-algebra and that we are given a system of roots of $q$ indexed by $D$ in $R$.
Then, there exists a unique morphism of $K$-algebras 
\begin{displaymath}
\xymatrix@R0cm{
K[t^{\frac 1D}] \ar[r] & R
\\t^r \ar@{|->}[r] & q^r.
}
\end{displaymath}
When $q \in R^\times$, it extends uniquely to $K[t^{\pm\frac 1D}]$.
\end{prop}

\begin{pf}
The morphism of monoids $N \to R$ will extend uniquely to a morphism of $K$-algebras.
The same result holds with $\pm N$ when $q \in R^\times$.
$\quad \Box$
\end{pf}

\begin{dfn} \label{defad}
A system $\{q_{n}\}_{n \in D}$ of roots of $q$ is said to be \emph{admissible} if
\begin{displaymath}
\forall n \in D, \quad (n)_{q_n} \in R^\times.
\end{displaymath}
\end{dfn}

\begin{xmp}
\begin{enumerate}
\item If $1 - q \in R^\times$, then any system of roots of $q$ is admissible because
\begin{displaymath}
(1 -q_{n})(n)_{q_{n}} = 1 -q_{n}^n = 1 -q.
\end{displaymath}
This applies in particular when $R$ is a field and $q \neq 1$, or more generally when $R$ contains a field $K$ and $q \in K$ with $q \neq 1$.
\item A non trivial system of roots of $1$ \emph{cannot} be admissible: we will have $(n)_{q_n} = 0$ for all $n \in D$.
\item Assume $R = \mathbb Z[\sqrt{-1}]$ and $q = -1$. Then the square roots of $q$ are not admissible because $1 \pm \sqrt{-1}$ is not invertible in $R$.
\end{enumerate}
\end{xmp}

\begin{dfn}
Let $D \subset \mathbb N \setminus \{0\}$ be a full inductive set of denominators for a submonoid $N \subset \mathbb Q_{\geq 0}$.
Let $\{q_{n}\}_{n \in D}$ be an admissible system of roots of $q$ in $R$.
If $r = \frac mn \in N$ with $m \in \mathbb N$ and $n \in D$, then the \emph{$q$-state} (or \emph{quantum state}) of $r$ is
\begin{displaymath}
(r)_{q} := \frac { (m)_{q_{n}} } { (n)_{q_{n}} }.
\end{displaymath}
If $q \in R^\times$ and $r \in \pm N$, then the \emph{$q$-state} (or \emph{quantum state} ) of $r$ is defined by the same formula.
\end{dfn}

Note that we must verify that this definition makes sense.
More precisely, since $D$ is assumed to be inductive, we must check that if $k \in \mathbb N$ is such that $kn \in D$, we also have
\begin{displaymath}
(r)_{q} = \frac { (km)_{q_{kn}} } { (kn)_{q_{kn}} }.
\end{displaymath}
But we know from proposition \ref{arit} that
\begin{displaymath}
(km)_{q_{kn}} = (k)_{q_{kn}}(m)_{q_{n}} \quad \mathrm{and} \quad  (kn)_{q_{kn}} = (k)_{q_{kn}}(n)_{q_{n}}.
\end{displaymath}

Following proposition \ref{monomap}, we will sometimes write $q^{\frac mn} := q_{n}^m$.
Then, the formula reads 
\begin{displaymath}
\left(\frac mn\right)_q := \frac { \sum_{i=0}^{m-1} q^{\frac in} } { \sum_{i=0}^{n-1} q^{\frac in} },
\end{displaymath}
and for example, we will have
\begin{displaymath}
\left(\frac 12\right)_{q} = \frac 1{1 + q^{\frac 12}}, \quad \left(\frac 13\right)_{q} = \frac 1{1 + q^{\frac 13} + q^{\frac 23}}, \quad \left(\frac 23\right)_{q} = \frac {1 + q^{\frac 13}} {1 + q^{\frac 13} + q^{\frac 23}}, \quad \ldots
\end{displaymath}
Also, if $q \in R^\times$, we will have
\begin{displaymath}
\left(-\frac mn\right)_{q} := -\frac { \sum_{i=1}^{m} q^{-\frac {i}n} } { \sum_{i=0}^{n-1} q^{\frac in} },
\end{displaymath}
and in particular
\begin{displaymath}
\left(-\frac 12\right)_{q} = -\frac 1{q^{\frac 12}+q}, \quad \left(-\frac 13\right)_{q} = -\frac 1{ q^{\frac 13} + q^{\frac 23}+q}, \quad \ldots
\end{displaymath}

Propositions \ref{arit} and \ref{pqzer} generalize as follows:

\begin{prop} \label{aritn}
Assume that we are given an admissible system of roots of $q$ indexed by an inductive set of denominators $D$ of a submonoid $N$ of $\mathbb Q_{\geq 0}$.
Then, for all $r \in N$ (or $r \in \pm N$ when $q \in R^\times$), we have
\begin{displaymath}
(1 -q)(r)_{q} = 1 - q^r.
\end{displaymath}
In particular, if $1 - q$ is invertible in $R$, we have
\begin{displaymath}
(r)_{q} = \frac {1 - q^r}{1 -q}.
\end{displaymath}
\end{prop}

\begin{pf}
If we write $r = \frac mn$ with $m \in N$ and $n \in D$, we have
\begin{displaymath}
(1 -q)(r)_{q} = \left( (1 -q^{\frac 1n}) (n)_{q^{\frac 1n}}\right)\left( \frac { (m)_{q^{\frac 1n}} } { (n)_{q^{\frac 1n}} }\right)
\end{displaymath}
\begin{displaymath}
= (1 -q^{\frac 1n})  (m)_{q^{\frac 1n}} = 1 - (q^{\frac 1n})^m = 1 - q^r. \quad \Box
\end{displaymath}
\end{pf}

\begin{prop}
Assume that we are given an admissible system of roots of $q$ indexed by an inductive set of denominators of a submonoid $N$ of $\mathbb Q_{\geq 0}$.
For all $r_1, r_2 \in N$ (or $\pm N$ when $q \in R^\times$), we have
\begin{displaymath}
(r_1 + r_2)_{q} = (r_1)_{q} + q^{r_1}(r_2)_{q}
\end{displaymath}
and
\begin{displaymath}
(r_1r_2)_{q} = (r_1)_{q} (r_2)_{q^{r_1}}.
\end{displaymath}
\end{prop}

\begin{pf}
We easily reduce to the case $R = \mathbb Q(t^{\frac 1D})$ and $q = t$ in which case, proposition \ref{aritn} tells us that
\begin{displaymath}
(r)_{q} = \frac {1 - q^r}{1 -q}
\end{displaymath}
whenever $r \in N$ (or $\pm N$ when $q \in R^\times$).
Then, the formulas are easily checked exactly as in the the proof of proposition \ref{pqzer} (integer case).
$\quad \Box$
\end{pf}

\begin{rmk} (see definition 2.1 of \cite{Ernst06} for example)
If $q \in \mathbb R_{>0}$ is not equal to $1$, one defines the \emph{$q$-analog} (or \emph{quantum analog}) of a real number $a$ as
\begin{displaymath}
(a)_{q} = \frac {1-q^a}{1-q}
\end{displaymath}
This is compatible with the above definition of the $q$-state of $r$ when $a = r \in \mathbb Q$.
More generally, if we choose a branch of the logarithm which is defined at a complex number $q \neq 1$, we may define the quantum analog of a complex number $a$ with the same formula and the convention $q^a = \exp(a\ln(q))$.
There are analogous results over ultrametric fields.
\end{rmk}

\section{Twisted powers} \label{twisted}

We assume here that $R$ is commutative and we fix a commutative $R$-algebra $A$ endowed with an $R$-algebra endomorphism $\sigma$.

\begin{dfn} \label{twdef}
If $x \in A$ and $n \in \mathbb N$, the \emph{$n$-th twisted power} of $x$ (with respect to $\sigma$) is
\begin{displaymath}
x^{(n)_{\sigma}} := \prod_{i=0}^{n-1}\sigma^{i}(x).
\end{displaymath}
\end{dfn}

In other words, we have
\begin{displaymath}
x^{(0)_{\sigma}} = 1,\quad  x^{(1)_{\sigma}} = x,\quad  x^{(1)_{\sigma}} = x\sigma(x),\quad \ldots,\quad x^{(n)_{\sigma}} = x\sigma(x) \cdots \sigma^{n-1}(x),\quad \ldots
\end{displaymath}
The twisted powers can also be defined inductively by
\begin{displaymath}
x^{(n+1)_{\sigma}} = x^{(n)_{\sigma}} \sigma^{n}(x) = x\sigma(x^{(n)_{\sigma}}).
\end{displaymath}

\begin{xmp}
\begin{enumerate}
\item When $\sigma(x) = x$, we have $x^{(n)_{\sigma}} = x^n$.
In particular, in the case $\sigma = \mathrm{Id}_{A}$, twisted powers are just usual powers.
\item
\begin{enumerate}
\item If $\sigma(x) = x - 1$, we obtain the \emph{falling Pochhammer symbol} of $x$:
\begin{displaymath}
x^{(n)_{\sigma}} = x(x-1) \cdots (x-n+1).
\end{displaymath}
When $R$ is a $\mathbb Q$-algebra, it is then common to extend the usual binomial coefficients by writing
\begin{displaymath}
{x \choose n} := \frac {x^{(n)_{\sigma}} }{n!}.
\end{displaymath}
\item If $\sigma(x) = x + 1$, we obtain the \emph{rising Pochhammer symbol} of $x$:
\begin{displaymath}
x^{(n)_{\sigma}} = x(x+1) \cdots (x+n-1)
\end{displaymath}
and:
\begin{displaymath}
\frac {x^{(n)_{\sigma}} }{n!} =  {x+n-1 \choose n}
\end{displaymath}
when $R$ is a $\mathbb Q$-algebra.
\end{enumerate}
\item More generally, if $\sigma(x) = x - h$ (resp. $\sigma(x) = x + h$) with $h \in R^\times$ and $R$ is a $\mathbb Q$-algebra, we will have the identity
\begin{displaymath}
\frac {x^{(n)_{\sigma}} }{n!} = h^n {x/h \choose n} \quad (\mathrm{resp.} \quad \frac {x^{(n)_{\sigma}} }{n!} = h^n {(x+n-1)/h \choose n}).
\end{displaymath}
\item Assume now that $\mathrm{Char}(R) = p > 0$ and $\sigma(x) = x + h$ with $h \in R$.
Then
\begin{displaymath}
x^{(p)_{\sigma}} = x^p -h^{p-1}x.
\end{displaymath}
When $h = 1$, this is exactly the \emph{Artin-Schreier} map.
\item If $\sigma(x) = qx$ with $q \in R$, then
\begin{displaymath}
(1-x)^{(p)_{\sigma}} = (1-x)(1-qx)\cdots (1-q^{p-1}x)
\end{displaymath}
is the \emph{$q$-Pochhammer symbol} of $x$.
\item If $y \in A$ satisfies $\sigma(y) = qy$ with $q \in R$, we may endow the polynomial ring $A[\xi]$ with the endomorphism $\sigma(\xi) = \xi + y$.
Then, we will have
\begin{displaymath}
\forall n \in \mathbb N, \quad \xi^{(n)_{\sigma}} = \xi(\xi + y) \cdots (\xi + (n-1)_{q}y).
\end{displaymath}
These twisted powers play an important role in the theory of $q$-difference equations.
\end{enumerate}
\end{xmp}

\begin{prop}
Assume $\sigma(x) = qx$ with $q \in R$.
Then, we have the following:
\begin{enumerate}
\item $x^{(n)_{\sigma}} = q^{\frac {n(n-1)}2}x^n$.
\item If $q\mathrm{-char}(R) = p$ is an \emph{odd} integer, then $x^{(p)_{\sigma}} = x^p$.
\item If $R$ is $q$-flat and $q\mathrm{-char}(R) = p > 0$, then  $x^{(p)_{\sigma}} = (-1)^{p-1}x^p$.
\end{enumerate}
\end{prop}

\begin{pf}
We have for all $i \in \mathbb N$, $\sigma^i(x) = q^ix$ and therefore
\begin{displaymath}
x^{(n)_{\sigma}} = \prod_{i=0}^{n-1} q^ix = q^{\frac {n(n-1)}2}x^n.
\end{displaymath}

If $q\mathrm{-char}(R) = p$, we have $q^p = 1$.
Therefore, if $p = 2k+1$, we obtain
\begin{displaymath}
x^{(p)_{\sigma}} = q^{\frac {p(p-1)}2} x^p= q^{kp} x^p= x^p.
\end{displaymath}

Finally, assume that $R$ is $q$-flat and $q\mathrm{-char}(R) = p > 0$.
For $p$ odd, we just proved the formula.
On the other hand, if $p = 2k$, we know from proposition \ref{even} that $q^{k} = -1$.
It follows that
\begin{displaymath}
x^{(p)_{\sigma}} = q^{k(p-1)} x^p= (-1)^{p-1} x^p. \quad \Box
\end{displaymath}
\end{pf}

\begin{lem} \label{sigit}
We have
\begin{enumerate}
\item $\forall x \in A, \forall n, m \in \mathbb N, \quad x^{(n)_{\sigma}}\sigma^n(x^{(m)_{\sigma}}) =  x^{(n+m)_{\sigma}}$
\item $\forall x, y \in A, \forall n \in \mathbb N \quad (xy)^{(n)_{\sigma}} = x^{(n)_{\sigma}}y^{(n)_{\sigma}}$
\item $\forall x \in A, \forall n,k \in \mathbb N \quad \sigma^k(x^{(n)_{\sigma}}) = \sigma^k(x)^{(n)_{\sigma}}$
\end{enumerate}
\end{lem}

\begin{pf}
All the equalities follow from the fact that $\sigma$ is a ring endomorphism.
More precisely, we have
\begin{displaymath}
x^{(n)_{\sigma}}\sigma^n(x^{(m)_{\sigma}}) =  \prod_{i=0}^{n-1}\sigma^{i}(x) \sigma^n\left(\prod_{j=0}^{m-1}\sigma^{j}(x)\right)
\end{displaymath}
\begin{displaymath}
=  \prod_{i=0}^{n-1}\sigma^{i}(x)\prod_{j=0}^{m-1}\sigma^{n+j}(x) = \prod_{i=0}^{m+n-1}\sigma^{i}(x) = x^{(n+m)_{\sigma}}.
\end{displaymath}
Also,
\begin{displaymath}
(xy)^{(n)_{\sigma}} = \prod_{i=0}^{n-1}\sigma^{i}(xy) =  \prod_{i=0}^{n-1}\sigma^{i}(x)\sigma^i(y) = \prod_{i=0}^{n-1}\sigma^{i}(x)\prod_{i=0}^{n-1}\sigma^{i}(y) = x^{(n)_{\sigma}}y^{(n)_{\sigma}}.
\end{displaymath}
And finally,
\begin{displaymath}
\sigma^k(x^{(n)_{\sigma}}) = \sigma^k(\prod_{i=0}^{n-1}\sigma^{i}(x))  = \prod_{i=0}^{n-1}\sigma^{i+k}(x)  = \sigma^k(x)^{(n)_{\sigma}}. \quad \Box
\end{displaymath}
\end{pf}

There is also a formula for moving from $\sigma$ to $\sigma^m$ that is quite useful:

\begin{prop} \label{mov}
We have
\begin{displaymath}
\forall x \in A, n, m \in \mathbb N, \quad \left( x^{(n)_{\sigma^m}}\right)^{(m)_{\sigma}} = \left( x^{(n)_{\sigma}}\right)^{(m)_{\sigma^n}} = x^{(mn)_{\sigma}}
\end{displaymath}
\end{prop}

\begin{pf}
We simply do the computations.
We have
\begin{displaymath}
\left( x^{(n)_{\sigma^m}}\right)^{(m)_{\sigma}} = \prod_{i=0}^{m-1}\sigma^{i} \left(\prod_{j=0}^{n-1}(\sigma^m)^j(x)\right)
\end{displaymath}
\begin{displaymath}
 =  \prod_{i=0}^{m-1} \prod_{j=0}^{n-1}\sigma^{mj+i}(x) =  \prod_{k=0}^{mn-1} \sigma^{k}(x) = x^{(mn)_{\sigma}}.
\end{displaymath}
And
\begin{displaymath}
\left( x^{(n)_{\sigma}}\right)^{(m)_{\sigma^n}} = \prod_{i=0}^{m-1}(\sigma^n)^i \left(\prod_{j=0}^{n-1}\sigma^j(x)\right)
\end{displaymath}
\begin{displaymath}
 =  \prod_{i=0}^{m-1} \prod_{j=0}^{n-1}\sigma^{ni+j}(x) =  \prod_{k=0}^{mn-1} \sigma^{k}(x) = x^{(mn)_{\sigma}}. \quad \Box
\end{displaymath}
\end{pf}

The \emph{twisted binomial formula} reads as follows:

\begin{prop} \label{twisbin}
Assume that $x, y \in A$ satisfy $\sigma(x) = qx$ with $q \in R$, and $\sigma(y) = y$.
Then, we have
\begin{displaymath}
\forall n \in \mathbb N, \quad (x+y)^{(n)_{\sigma}} = \sum_{k=0}^n {n \choose k}_q x^{(k)_{\sigma}}y^{(n-k)_{\sigma}}.
\end{displaymath}
\end{prop}

\begin{pf}
This is exactly the formula of proposition \ref{bincof}.
$\quad \Box$
\end{pf}

We can also state the \emph{Frobenius property}:

\begin{cor}
Let $q \in R$.
Assume $q\mathrm{-char}(R) = p > 0$ and $R$ is $q$-flat.
Assume that $x, y \in A$ satisfy $\sigma(x) = qx$ and $\sigma(y) = y$.
Then, we have
\begin{displaymath}
\forall n \in \mathbb N, \quad (x+y)^{(p)_{\sigma}} = x^{(p)_{\sigma}} + y^{(p)_{\sigma}}
\end{displaymath}
\end{cor}

\begin{pf}
Follows from lemma \ref{qbin}.
$\quad \Box$
\end{pf}

Note that one can also recover the quantum Lucas theorem (proposition \ref{binexp}) as a corollary of proposition \ref{twisbin} as explained for example in lemma 1 of \cite{littelmann98}.

We also want to mention that one can consider the notion of twisted powers of an ideal (products and images are always meant as ideals):

\begin{dfn} \label{twdefid}
If $\mathfrak a \subset A$ is an ideal, the \emph{twisted powers} of $\mathfrak a$ are defined as
\begin{displaymath}
\mathfrak a^{(n)_{\sigma}} := \prod_{i=0}^{n-1}\sigma^{i}(\mathfrak a).
\end{displaymath}
And the \emph{twisted completion} of $A$ along $\mathfrak a$ is
\begin{displaymath}
\hat A^\sigma := \varprojlim A/\mathfrak a^{(n)_{\sigma}} 
\end{displaymath}
\end{dfn}

Again, it means that
\begin{displaymath}
\mathfrak a^{(0)_{\sigma}} := A, \quad \mathfrak a^{(1)_{\sigma}} := \mathfrak a, \quad \mathrm{and} \quad   \mathfrak a^{(n)_{\sigma}} = \mathfrak a\sigma(\mathfrak a) \cdots \sigma^{n-1}(\mathfrak a).
\end{displaymath}

Of course, when $\mathfrak a$ is a principal ideal, say $\mathfrak a = (x)$, we will have $\mathfrak a^{(n)_{\sigma}} = (x^{(n)_{\sigma}})$.

\begin{xmp}
\begin{enumerate}
\item If $\sigma = \mathrm{Id}_{A}$ these are just usual powers and usual completion.
\item If $R = \mathbb Q$, $A = \mathbb Q [x]$ and $\sigma(x) = x + h$ with $h \neq 0$, then $\hat A^\sigma \simeq \mathbb Q[x]^{\mathbb N}$ and the canonical map $A \to \hat A^\sigma$ sends $x$ to $(x-ih)_{i \in \mathbb N}$.
\end{enumerate}
\end{xmp}

\bibliographystyle{plain}
\addcontentsline{toc}{section}{References}
\bibliography{BiblioBLS}

\end{document}